\documentclass[12pt]{article}
\usepackage{amsfonts}
\usepackage{amsmath, graphicx, epsfig,psfrag, verbatim}
\usepackage{eucal,  amssymb,amsmath,amsthm,graphicx, amsfonts, latexsym}
\usepackage[utf8]{inputenc}
\newtheorem{theorem}{Theorem}[section]

\newtheorem{cor}{Corollary}[section]

\newtheorem{defin}{Definition}[section]
\newcommand{\floor}[1]{\lfloor #1 \rfloor}

\begin{document} \parskip=5pt plus1pt minus1pt \parindent=0pt
\title{Directed preferential attachment models}
\author{Tom Britton$^1$, Stockholm University}
\date{\today}
\maketitle

\begin{abstract}
The directed preferential attachment model is revisited. A new exact characterization of the limiting in- and out-degree distribution is given by two \emph{independent} pure birth processes that are observed at a common exponentially distributed time $T$ (thus creating dependence between in- and out-degree). The characterization gives an explicit form for the joint degree distribution, and this confirms previously derived tail probabilities for the two marginal degree distributions. The new characterization is also used to obtain an explicit expression for tail probabilities in which both degrees are large. A new generalised directed prefererantial attachment model is then defined and analysed using similar methods. The two extensions, motivated by empirical evidence, are to allow double-directed (i.e.\ undirected) edges in the network, and to allow the probability to connect an ingoing (outgoing) edge to a specified node to also depend on the out-degree (in-degree) of that node.

\end{abstract}


\footnotetext[1]{Department of Mathematics, Stockholm University, 106 91 Stockholm, Sweden.\\ Email: tom.britton@math.su.se}

\section{Introduction and models}\label{sec-intro}

The (undirected) preferential attachment model (PA) is a random network model defined by Barab{\'a}si and Albert \cite{BA99}. To start off, the network consists of one single node without any edge. At each time step $k=1,2,\dots ,$ a new node with $m$ (a fixed integer-parameter in the model) new edges connected to it, is added. Each of the new edges of the node is connected, independently, to existing nodes, and the probability to connect to a specific node with current degree $i$ is proportional to $i$. Two novel features, as compared to most other network models at the time, were that it was defined sequentially, thus with nodes having different ages, and that the degree distribution of nodes in a large network were shown to have power law tails rather than exponentially decaying tail probabilities.

In 2003, Bollob{\'a}s et al \cite{B03} defined a related model, but now for a network in which edges are directed rather than undirected. 
As in the undirected model, edges/nodes are entered at each discrete time step. However, now one of three different possibilities may happen: 1) either a new node with an \emph{outgoing} edge is added (with probability $\alpha$), or 2) a new directed \emph{edge} but no new node is added (with probability $\beta$), or 3) a node with an \emph{ingoing} edge is added (with probability $\gamma=1-\alpha-\beta$). In the first event, the edge connects to (i.e.\ points at) a node $u$ having current in- and out-degree $i$ and $j$ with probability proportional to $i+\delta_I$. In the second event the edge starts from $u$ with probability proportional to $j+\delta_O$ and, independently, ends at $u$ with probability proportional to $i+\delta_I$. In the third and final event the edge starts from $u$ with probability proportional to $j+\delta_O$. The starting configuration is not important for limiting properties of the network as long as there are few edges and nodes, but here we define it by a single node without any edges.

The model has 4 parameters: $\alpha,\ \beta, \ \delta_I,\ \delta_O$ (recall that $\gamma=1-\alpha-\beta$). It is important that $\delta_I>0$ and $\delta_O>0$; otherwise nodes born with in-degree 0 will never get positive in-degree, and similarly for the remaining nodes born with out-degree 0. In order to avoid some less interesting special cases of the model (which would require special attention in the analyses) we will also assume that $\alpha >0$ and $\gamma >0$. We denote this model the Directed Preferential Attachment (DPA) model.

In Bollob{\'a}s et al \cite{B03} it was shown that $E(X_i(n))=\varphi_in+o(n)$, $E(Y_j(n))=\phi_jn+o(n)$ and $E(M_{ij}(n))=f_{ij}n + o(n)$, where $ X_i(n)$ denotes the number of nodes with in-degree $i$, $Y_j(n)$ denotes the number of nodes with out-degree $j$, and $M_{ij}(n)$ denotes the number of nodes having in-degree $i$ and out-degree $j$, at time step $n$. Bollob{\'a}s et al \cite{B03} did not obtain explicit expressions for $\varphi_i$, $\phi_j$ or $f_{ij}$, but instead derived limiting properties for large $i$ and $j$ (sending either, but not both, off to infinity for $f_{ij}$). More specifically, they show that $\varphi_i\sim i^{-(1+1/c_I)}$ and $\phi_j\sim j^{-(1+1/c_O)}$, where
\begin{align}
c_I &=\frac{\alpha+\beta}{1+\delta_I(\alpha+\gamma)}\label{c_I}
\\
c_O &= \frac{\gamma+\beta}{1+\delta_O(\alpha+\gamma)},\label{c_O}
\end{align}
 ($a_i\sim b_i$ meaning $a_i/b_i\to c$ for some $0<c<\infty$ as $i\to\infty$). As for the bivariate distribution $f_{ij}$ it was shown that, having one of $i$ and $j$ fixed and letting the other tend to infinity, they showed that $f_{ij}\sim i^{-x_I}$ (as $i\to\infty$ and $j$ fixed) and $f_{ij}\sim j^{-x_O}$ (as $j\to\infty$ and $i$ fixed), where $x_I= 1+1/c_I+c_O(\delta_O+I_{(\gamma\delta_O=0)})/c_I$, and $x_O= 1+1/c_O+c_I(\delta_I+I_{(\alpha\delta_I=0)})/c_O$ (note that there is a small misprint the first time the expressions appear in Bollob{\'a}s et al, \cite{B03}). The method used when proving the results is by analysing the evolution of the vector valued processes $(X_0(n), X_1(n), \dots ),\ (Y_0(n), Y_1(n),\dots )$ and $(M_{10}(n), M_{01}(n), \dots )$, and by analysing the limiting partial differential equations. Samorodnitsky et al.\ (2016) take this one step further and derive an exact integral characterisation of the joint generating function $\varphi (x,y)$ of $\{ f_{ij}\}$. Using this characterization they are also able to prove that the joint distribution $\{ f_{ij}\}$ has jointly regularly varying tails with a specified tail measure.

In the current paper we analyse the same DPA model but using a different method. Instead we analyse the evolution of the in- and out-degree of \emph{one} randomly selected node (born before $n$) up until time step $n$. If we let $p_{ij}^{(n)}$ denote the probability that the degree of this node equals $(i,j)$ at time step $n$, it follows that $\lim_n p_{ij}^{(n)} =cf_{ij}$ (the constant $c$ is only there because the $f_{ij}$ of Bollob{\'a}s et al. \cite{B03} will not sum to unity). Using this alternative method we show, by  means of weak convergence and a certain time transformation, that the evolution of degrees of a randomly selected node, in the limit, converges to the evolution of two independently evolving Markov birth processes, which are both stopped at a common expontially distributed time $T$. We now define this limiting process.

Consider a bivariate birth process $X(t)=(X_I(t),X_O(t))$ with time-homogeneous birth rates 
\begin{align*}
P(X(t+h)=(i+1,j)|X(t)=(i,j)) &=\lambda_{ij}^Ih +o(h)
\\
P(X(t+h)=(i,j+1)|X(t)=(i,j)) &=\lambda_{ij}^Oh +o(h)
\\
P(X(t+h)=(i,j)|X(t)=(i,j))&=1-\lambda_{ij}^Ih - \lambda_{ij}^Oh +o(h),
\end{align*}
where the jump intensities are given by
\begin{align}
\lambda_{ij}^I &= \lambda_{i}^I = (i+\delta_I)c_I\label{lambda_i}
\\
\lambda_{ij}^O &= \lambda_{j}^O =(j+\delta_O)c_O,\label{lambda_j}
\end{align}
with $c_I$ and $c_O$ defined in (\ref{c_I}) and (\ref{c_O}).

Because the two jump rates depend only on the first and second coordinate respectively, the two coordinate processes, reflecting in- and out-degree, evolve \emph{independently}. Assume that the process is started either in state (0,1) or (1,0): 
\begin{align*}
P(X(0)=(0,1))&= \frac{\alpha}{\alpha+\gamma},\text{ and} \\
P(X(0)=(1,0))&= \frac{\gamma}{\alpha+\gamma}. 
\end{align*}

In Section \ref{sec-results} we make use of this process and shown that the limiting degree distribution of the DPA model is identical to that of this process observed after an exponentially distributed time, which enables the computation of its bivariate tail distribution. Another advantage with this new characterization is that it easily extends to related, but more realistic preferential attachment models, including the one defined below.

The DPA-model of Bollob{\'a}s et al.\ \cite{B03} has two main features which may be critizised from the point of realism. The first is that the two different degrees of a node evolve independently in the sense that the rate/probability of aquiring an additional ingoing edge (thus increasing the in-degree by 1) depends on the current in-degree of the node, but not on its out-degree, and vice versa. A more realistic model would in many situations be to let the probability that an added directed edge points to a node having current in- and out-degree $(i,j)$ to be proportional to $i+cj+\delta_I$, and similarly that the probability that a new edge points out from this node to be proportional to $di+j+\delta_O$ (where $c$ and $d$ are non-negative model parameters; in the original DPA-model $c=d=0$). This will allow for a stronger dependence between in- and out-degree. 

A second, perhaps more important, feature that the original DPA-model can be criticized for is that the fraction of edges that are ''double directed'' (or equivalently undirected), i.e.\ pair of nodes for which there exist directed edges going both ways between them, will be negligible. In many empirical networks having directed edges, the fraction $\rho$ of directed edges for which the reciprocal edge is also present, is far away from 0, typically in the interval $(0.2,\ 0.8)$ (cf.\ Stanford data base \cite{LK14}, and Spricer and Britton \cite{SB15} who consider another model for a partially directed random network). This can be achieved if we modify the DPA model by simply stating that, at each time step (when a directed edge is added), the corresponding reciprocal edge is added with probability $\rho$. 

From the reasoning above we now define what we call the Generalised Directed Preferential Attachment (GDPA) model.

\begin{defin}[The Generalised Directed Preferential Attachment (GDPA) model] \label{def-gdpa}
The process is started at $k=0$ with a single node without any edges. At each discrete time step $k=1,2,\dots$, one of three different events can happen: 1) with probability $\alpha$ a new node with an edge pointing out from this node is added, 2) with probability $\beta$ a new directed edge without nodes is added, and 3) with the remaining probability $\gamma=1-\alpha-\beta$ a new node with a directed edge pointing at the new node is added. In the first and second cases, the probability that the new edge points to a specific existing node with in- and out-degree $(i,j)$ is proportional $i+cj+\delta_I$. In the second and third case, the probability that the new edge points out from a specific node with in- and out-degree $(i,j)$ is proportional $di+j+\delta_O$. Finally, in every step, the added directed edge is made reciprocal (i.e.\ double-directed) with probability $\rho$, independently in each time step. 
\end{defin}
  
The GDPA model has 7 parameters: $\alpha,\ \beta, \ \delta_I,\ \delta_0, c,\ d$ and $\rho$ (recall that $\gamma=1-\alpha-\beta$ and hence not a free parameter). Note that in the DPA model it was crucial that $\delta_I>0$ and $ \delta_O>0$, since otherwise a node would only have positive in- or out-degree. For the GDPA model this is no longer necessary when $c>0$ and $d>0$. For this reason it is possible to reduce the number of parameters to 5 by assuming $\delta_I= \delta_O=0$. Further, the GDPA model is identical to the DPA model when instead $c=d=\rho=0$.

\section{Main results}\label{sec-results}

We now state our main result characterising the limiting degree distribution of the directed preferential attachment model in terms of our simple 2-dimensional birth process $X(\cdot)$ (defined in the previous section) evaluated after and exponential time. Recall that we assume $\alpha,\ \gamma,\ \delta_I$ and $\delta_O$ to all be strictly positive (the interesting case) to avoid special cases.
\begin{theorem}
\label{mainth}
Let $p_{ij}^{(n)}$ denote the probability that a randomly selected node in the DPA model after $n$ steps has in- and out-degree $i$ and $j$ respectively. Let $X(t)$ denote the bivariate birth process defined above, and $T$ an independent Exp(1) random variable. Then,
$$
p_{ij}^{(n)}\to p_{ij}:=P(X(T)=(i,j)).
$$
\end{theorem}

\textbf{Remark}. Let $N(n)$ denote the number of nodes after $n$ steps in the DPA model (which follows $1+Bin(n,\alpha+\gamma)$ distribution). It then directly follows that $M_{ij}(n)/N(n)$, the fraction of nodes having in- and out-degree $(i,j)$, converges in probability to $p_{ij}$.
\vskip.3cm

The explicit distribution of $X(T)$ can be derived in two ways. Either by studying the embedded discrete-time random walk which, if currently in state $(k,\ell)$, goes to state $k+1,\ell$ with probability $\lambda_k^I/(\lambda_k^I+\lambda^O_{\ell}+1)$, to state $k,\ell+1$ with probability $\lambda_{\ell}^O/(\lambda_k^I+\lambda^O_{\ell}+1)$ or gets stuck for ever in state $(k,\ell)$ with the remaining probability $1/(\lambda_k^I+\lambda^O_{\ell}+1)$. The probability $P(X(T)=(i,j))$ is then the probability that this random walk gets stuck in state $(i,j)$. Note that there are many different paths, in general having different probabilities, ending in state $(i,j)$. This derivation technique will be applied when analysing the GDPA model.

The other way to derive the distribution is by integrating over possible values of $T$, and using that the two components of the continuous time bivariate processes evolve independently. Let $X^{(0,1)}(t)=(X_I^0(t),X_O^1(t))$ and $X^{(1,0)}(t)=(X_I^1(t),X_O^0(t))$ denote the bivariate birth process described above, but where we condition on the starting state being $(0,1)$ and $(1,0)$ respectively (remember that the probabilities for these two starting points are $\alpha/(\alpha +\gamma)$ and $\gamma/(\alpha +\gamma)$ respectively). We then have
\begin{align}
p_{ij} & = \frac{\alpha}{\alpha+\gamma} \int_0^\infty P(X_I^0(t)=i)P(X_O^1(t)=j)e^{-t}dt \nonumber
\\
&\hskip3cm + \frac{\gamma}{\alpha+\gamma} \int_0^\infty P(X_I^1(t)=i)P(X_O^0(t)=j)e^{-t}dt,\label{p_ijexact}
\end{align}
and the following result gives explicit expressions for the bivariate degree distribution (using the notation $a_k=a(a-1)\dots (a-k+1)$).

\begin{cor}\label{cor-lim-deg-dist}
The marginal distributions of the two birth processes conditioned on their starting values $r$, are given by
\begin{align*}
P(X_I^r(t)=i) &= \frac{(\delta_I+i-1)_{i-r}}{(i-r)!}  e^{-c_I(\delta_I+r)t}\left(1-e^{-c_It}\right)^{i-r} ,\ i=r,r+1,\dots
\\
P(X_O^r(t)=j) &= \frac{(\delta_O+j-1)_{j-r}}{(j-r)!}  e^{-c_O(\delta_O+r)t}\left(1-e^{-c_Ot}\right)^{j-r},\ j=r, r+1, \dots .
\end{align*}
The marginal distributions of the stopped birth processes conditioned on their starting value $r$, are given by
\begin{align*}
P(X_I^r(T)=i) &= \frac{(\delta_I+i-1)_{i-r}}{(\delta_I+\frac{1}{c_I} + i)_{i-r+1}}   
\\
P(X_O^r(T)=j) &= \frac{(\delta_O+j-1)_{j-r}}{ \delta_O+\frac{1}{c_O} + j)_{j-r+1} }  .
\end{align*}

The joint degree distribution $p_{ij}=P(X(T)=(i,j))$ is given by
\begin{align*}
p_{ij}&= \frac{\alpha}{\alpha+\gamma} \frac{(\delta_I+i-1)_i(\delta_O+j-1)_{j-1}}{i!(j-1)!}\sum_{k=0}^i\sum_{\ell=0}^{j-1} (-1)^{k+\ell }
\frac{\binom{i}{k}\binom{j-1}{\ell}}{c_I(\delta_I+k)+c_O(\delta_O+\ell)+1},
\\
 &\qquad + \frac{\gamma}{\alpha+\gamma}\frac{(\delta_I+i-1)_{i-1}(\delta_O+j-1)_{j}}{(i-1)!j!} \sum_{k=0}^{i-1}\sum_{\ell=0}^j (-1)^{k+\ell }
\frac{\binom{i-1}{k}\binom{j}{\ell}}{c_I(\delta_I+k)+c_O(\delta_O+\ell)+1} .
\end{align*}

\end{cor}

Beside having an explicit (albeit long) form for the limiting degree distribution $\{ p_{ij}\}$, its characterization as a stopped bivariate birth process making independent jumps component-wise, gives hope for a richer analysis of the tail behavior of the two degrees. For example, it is not hard to confirm the earlier stated results of Bollob{\'a}s et al. \cite{B03}.

\begin{cor}[Bollob{\'a}s et al. \cite{B03}]\label{cor-Bollob}
For large in-degrees $p_i:=\sum_jp_{ij} \sim i^{-(1+1/c_I)}$, for large out-degrees $q_j:=\sum_ip_{ij}\sim j^{-(1+1/c_O)}$. For fixed $j$ and large $i$: $p_{ij}\sim i^{-x_I}$, and for fixed $i$ and large $j$: $p_{ij}\sim j^{-x_O}$, where $x_I= 1+1/c_I+c_O\delta_O/c_I$, and $x_O= 1+1/c_O+c_I\delta_I/c_O$.
\end{cor}

\textbf{Remark}. The original theorem of Bollob{\'a}s et al.\ contain an additional term for $x_I$ and $x_O$, but these vanish when restricting the parameters $\delta_I,\ \delta_O,\ \alpha$ and $\gamma$ all being strictly positive as we have done (the interesting case).
\vskip.3cm

Our new result concerns the tail probability $p_{ij}$ when both $i\to\infty $ and $j\to\infty$. More specifically, let $i=n$ and $j=\floor{sn^r}$ for some $0<s,r<\infty$, where $\floor{x}$ is the integer part of $x$. We have the following theorem.

\begin{theorem}\label{Theor-tail}
The tail probabilities $p_{ij}$, for $i=n$ and $j=\floor{sn^r}$, satisfy
\begin{equation}
p_{n,\floor{sn^r}}\sim n^{\delta_I-1+r(\delta_O-1)-(c_I\delta_I+\delta_Oc_O+1)\max (1/c_I, r/c_O)} 
=
 \begin{cases}
  n^{-\left( 1+\delta_O(c_O/c_I-r) + r + 1/c_I \right) } &\mbox{if } r\le  c_O/c_I,
  \\
  n^{-\left( 1+\delta_I(rc_I/c_O -1) + r +1/c_O \right) } &\text{if } r\ge  c_O/c_I.
 \end{cases}
\end{equation}
The factor $s$ hence has no effect on the tail. The choice of $r$ which maximizes the tail probability is given by
\begin{equation}
\textrm{Argmax}_{ \{r\ge 0\} }p_{n,\floor{n^r}} \to 
\begin{cases}
0 & \mbox{if }\delta_O<1,
\\
c_O/c_I & \mbox{if }\delta_O>1.
\end{cases}
\end{equation}
If $\delta_O=1$ any value of $r$ between 0 and $c_O/c_I$ give the same asymptotic tail. 
\end{theorem}

\vskip.3cm
The limiting degree distribution of the original (undirected) PA model of Albert and Barab{\'a}si \cite{BA99} can also be derived using similar methodology. Here the limiting process is even simpler: a 1-dimensional pure birth process $Z(t)$ starting at $Z(0)=m$ and having birth rate $\lambda_i=i/2$ that is stopped at $T\sim Exp(1)$. The limiting distribution has been derived earlier using other methods.

\begin{cor}[Dorogovtsev et al. \cite{DMS00}]\label{cor-Bol-Rio}
Let $p_i^{(n)}$ denote the probability that the degree of a randomly selected individual after $n$ steps has degree $i$ in the (undirected) PA-model, and let $Z(t)$ and $T$ be defined as above. Then 
$$
p_i^{(n)}\to p_i:= P(Z(T)=i)=\frac{2m(m+1)}{i(i+1)(i+2)},\ i=m,m+1,\dots.
$$
\end{cor}

We now consider the generalized directed preferential attachment (GDPA) model.
Using very similar methods as when proving Theorem \ref{mainth}, it can be shown that also for the GDPA model, the limiting degree distribution can be characterized by a continuous-time bivariate ''birth process'' which is stopped and observed after an exponentially distributed random time $T$. We now describe the limiting process to which the degree distribution of the GDPA model converges. Let $Y(t)=(Y_I(t),Y_O(t))$ be a time-homogeneous bivariate Markov birth process, but where now simultaneous births of the two components are also possible. The three different birth rates $\beta_{ij}^{I}$, $\beta_{ij}^{O}$ and $\beta_{ij}^{I+O}$ are defined by
\begin{align}
\beta_{ij}^I &= (i+cj+\delta_I)(1-\rho)\frac{\alpha+\beta}{(1+\rho)(1+c)+\delta_I(\alpha+\gamma)}=:\ (i+cj+\delta_I)(1-\rho)g_I \nonumber
\\
\beta_{ij}^O &= (di+j+\delta_O)(1-\rho)\frac{\gamma+\beta}{(1+\rho)(1+d)+\delta_O(\alpha+\gamma)}=:\ (di+j+\delta_O)(1-\rho) g_O \label{GDPA-int}
\\
\beta_{ij}^{I+O} &= (i+cj+\delta_I)\rho g_I + (di+j+\delta_O)\rho g_O ,\nonumber
\end{align}
$g_I$ and $g_O$ being the respective ratio expressions.
With the third rate is meant that $P(Y(t+h)=(i+1,j+1)|Y(t)=(i,j))= \beta_{ij}^{I+O}h+o(h)$. The process is started in either of the states $(0,1), (1,0)$ or $(1,1)$ with respective probability $(1-\rho)\alpha/(\alpha+\gamma)$, $(1-\rho)\gamma/(\alpha+\gamma)$ and $\rho$. As before, let $T\sim Exp(1)$ be an independent exponentially distributed time. We then have the following theorem.

\begin{theorem}
\label{ThGDPA}
Let $\pi_{ij}^{(n)}$ denote the probability that a randomly selected node in the GDPA model after $n$ steps has in- and out-degree $i$ and $j$ respectively. Let $Y(t)$ denote the bivariate birth process defined above, and $T$ an independent $Exp(1)$ random variable. Then,
$$
\pi_{ij}^{(n)}\to \pi_{ij}:=P(Y(T)=(i,j)).
$$
\end{theorem}

\textbf{Remark}. Just as in Theorem \ref{mainth}, it follows that the fraction of nodes having in- and out-degree $(i,j)$ converges in probability to $\pi_{ij}$.

\vskip.3cm

The limiting probabilities $\{\pi_{ij}\}$ can be computed by summing over all paths, starting from either $(0,1), (1,0)$ or $(1,1)$, and ending in $(i,j)$, where each jump either increases the in-degree by 1, the out-degree by 1, both degrees by 1, or makes a complete stop. If currently in state $(k,\ell )$, the process jumps to $(k+1,\ell)$ with probability $\beta^I_{k,\ell }/\sigma_{k,\ell}$, to state $(k,\ell+1)$ with probability $\beta^O_{k,\ell}/\sigma_{k,\ell}$, to state  $(k+1,\ell+1)$ with probability $\beta^{I+O}_{k,\ell }/\sigma_{k,\ell}$, or makes a complete final stop with probability $1/\sigma_{k,\ell}$, where $\sigma_{k,\ell} =\beta^I_{k,\ell}+\beta^O_{k,\ell}+\beta^{O+I}_{k,\ell}+1$ (the different $\beta$'s were defined in Equations (\ref{GDPA-int})).

As an illustration, 
\begin{align*}
\pi_{11} & = \frac{(1-\rho)\alpha}{\alpha+\gamma} \times \frac{\beta^I_{0,1 }}{\beta^I_{0,1 }+ \beta^O_{0,1 }+ \beta^{I+O}_{0,1 }+1}  \times \frac{1}{\beta^I_{1,1 }+ \beta^O_{1,1 }+ \beta^{I+O}_{1,1 }+1}
\\
& + \frac{(1-\rho)\gamma}{\alpha+\gamma}  \times \frac{\beta^O_{1,0 }}{\beta^I_{1,0 }+ \beta^O_{1,0 }+ \beta^{I+O}_{1,0 }+1}  \times \frac{1}{\beta^I_{1,1 }+ \beta^O_{1,1 }+ \beta^{I+O}_{1,1 }+1}
\\
& + \rho  \times \frac{1}{\beta^I_{1,1 }+ \beta^O_{1,1 }+ \beta^{I+O}_{1,1 }+1}.
\end{align*}
The first factor in each row is the probability of starting in $(0,1)$, $(1,0)$ and $(1,1)$ respectively. For higher degrees there will be many more paths to sum over. Starting at $(0,1)$ and ending in $(1,2)$ can for example happen in three different ways, either first jumping to $(1,1)$ followed by a jump to $(1,2)$ or first jumping to $(0,2)$ and then to $(1,2)$, or jumping directly to $(1,2)$. 

As for the tail probabilities of the GDPA model $\{\pi_{ij}\}$, it should be possible to derive them using a similar analysis as for the tails of the DPA model (Corollary \ref{cor-Bollob} and Theorem \ref{Theor-tail}). However, the fact that the two components no longer evolve independently makes the analysis more involved and its tail behaviour remains open problem.

\section{Proofs}

\subsection{Proof of Theorem \ref{mainth}}

The proof consists of two parts. First we show that the evolution of the degrees of a randomly selected node up until step $n$, converges to the evolution of degrees of a continuous time bivariate birth process $D=(D_I, D_O)$, born at $U\sim U[0,1]$, having \emph{time-inhomogeneous} birth rates $\lambda_{ij}^I/t$ and $\lambda_{ij}^O/t$ respectively, and observed at time $t=1$. Then, using a time transformation, we show that this limiting distribution has the same distribution as the \emph{time-homogeneous} bivariate process $X=(X_I,X_O)$ of the theorem, having birth rates $\lambda_{ij}^I$ and $\lambda_{ij}^O$ respectively, born at time 0 and observed at time $T\sim Exp(1)$.

Recall that $p_{ij}^{(n)}$ is the probability that a randomly chosen node after $n$ time steps has in- and out-degree $(i,j)$. At the start there are no edges, and at each time step, one edge is added, so there are $n$ edges at this time. Further, at the start there is 1 node, and at each time step a new node is added with probability $\alpha+\gamma$, so the number of nodes at time step $n$, denoted $N(n)$, is $1+Bin(n,\alpha + \gamma)$. For large $n$ this is well approximated by $(\alpha+\gamma)n$ which is done below. 

Fix $s$, $0\le s\le 1$, and consider a node that entered the network at time step $\floor{sn}$ (the integer part of $sn$). For $s\le t\le 1$ we define the bivariate process $(D_I^{(n)}(t),D_O^{(n)}(t))$ as the in- and out-degree of that node at time step $\floor{tn}$. The starting value of our node may either be $(0,1)$ or $(1,0)$ depending on if it entered through the event 1 or event 3 (if event 2 happens no new node is added). The probabilities are hence given by
$$
P(D_I^{(n)}(s)=0, D_O^{(n)}(s)=1)=\frac{\alpha}{\alpha+\gamma}= 1- P(D_I^{(n)}(s)=1, D_O^{(n)}(s)=0).
$$
Assume that at time $t$ our process equals $(D_I^{(n)}(t),D_O^{(n)}(t))=(i,j)$, and let $\Delta t=1/n$. We first compute the probability/intensity that our process will increase its in-degree by 1. This happens if either event 1 happens and our node is selected as the node to point at, or that event 2 happens and our node is selected as the node to point at. The probability for this is hence
\begin{align*}
P& (D_I^{(n)}(t+\Delta t),D_O^{(n)}(t+\Delta t)=(i+1,j)|\ D_I^{(n)}(t),D_O^{(n)}(t)=(i,j))
\\
&= (\alpha + \beta)\frac{i+\delta_I}{\floor{tn}+\delta_I N(\floor{tn})} +o(1/n)
\\
& = (\alpha + \beta)\frac{i+\delta_I}{t(1+\delta_I (\alpha+\gamma))}\Delta t +o(\Delta t).
\end{align*} 
For the out-degree we obtain the similar result
\begin{align*}
P & (D_I^{(n)}(t+\Delta t),D_O^{(n)}(t+\Delta t)=(i,j+1)|\ D_I^{(n)}(t),D_O^{(n)}(t)=(i,j))
\\
&= (\beta+\gamma)\frac{j+\delta_O}{\floor{tn}+\delta_O N(\floor{tn})} + o(1/n)
\\
& = (\beta+\gamma )\frac{j+\delta_O}{t(1+\delta_O (\alpha+\gamma))}\Delta t + o(\Delta t).
\end{align*}
It can in fact also happen that our node increases both its in- and out-degree. This happens in case event 2 happens and our node is chosen both for start and end of the edge (thus creating a loop). This event should in fact also be excluded from the above probabilities. However, this happens with a probability proportional to $1/n^2=(\Delta t)^2$, so these events will not occur in the limit.
Together with the observation that our process is tight it follows from standard results (e.g. \cite{Bil99}) that the process $(D^{(n)}_I,D^{(n)}_O)$ will, in the limit as $n\to\infty$, converge to a continuous-time stochastic process $(D_I, D_O)$ on $[s, 1]$. The limiting process only makes increases by 1 unit, and never both coordinates at the same time: it is a bivariate time-inhomogeneous Markov birth process.
An important observation is that the birth rate for the in-degree only depends on the current value of the in-degree $i$ but not on the current out-degree $j$, and similarly for the out-degree, which means that the two degrees evolve independently. We can hence drop this dependence in our notation and write $\lambda^{I}_{i}$ and $\lambda^{O}_{j}$ as defined in (\ref{lambda_i}) and (\ref{lambda_j}).
It hence follows that our process $(D_I^{(n)}(t),D_O^{(n)}(t))$ converges to a continuous time, time-inhomogeneous bivariate Markov birth process, with birth intensities given by $\lambda^I_i/t$ and $\lambda^O_j$ respectively.
The process starts at time $s$ with degrees $(0,1)$ or $(1,0)$, with respective probability $\alpha/(\alpha+\gamma)$ and $\gamma/(\alpha+\gamma)$.

Finally, our randomly chosen node has the same probability $\alpha+\gamma$ of being born any time point between 1 and $n$ implying that the birth time $s$ of the limiting process is $U[0, 1]$. We have thus shown that the distribution $\{ p_{ij}^{(n)}\}$ of the in- and out-degree of a randomly selected node after step $n$ in the GPA model, converges (as $n\to\infty$) to the distribution of $(D_I(1),D_O(1))$, where this process is started at a uniformly distributed time $s$ with starting configuration as stated above, and where  evolution of the two degrees evolves according to the rates given above until time $t=1$.

The second part of the proof consists of showing that this distribution equals the one stated in the theorem. We do this by looking at the accumulated intensities. Since the jump rates of $(D_I,D_O)$ are of the form $\lambda_{ij}\times 1/t$ and the jump rates of $(X_I,X_O)$ equal $\lambda_{ij}$ for the same $\lambda_{ij}$, it suffices to show that the accumulated rates agree for a fixed parameter $\lambda$ say.

The accumulated rate for the $D$ process, starting at time $s\sim U[0,1]$ point and evolving until time 1, equals
$$
\int_0^1 \left( \int_s^1 \lambda \frac{1}{t}dt\right)  ds =\lambda.
$$
The accumulated rate for the $X$ process starting at time 0 and evolving up until time $T\sim Exp(1)$ equals
$$
\int_0^\infty \left( \int_0^t\lambda ds\right) e^{-t}dt =\lambda .
$$
The two processes hence have the same accumulated jump rates and start from the same initial condition, which implies that they have the same distribution at the observation points. Another way to obtain this result is to make a time transformation of $D$: first we reverse time and let it evolve from 0 to $s\sim U[0,1]$ with rate $\lambda/(1-s)$, and then transform time: $s'=\log(1-s)$. These two changes combined lead to the rates and limits of the $X$-process. This completes the proof.

\subsection{Proof of Corollary \ref{cor-lim-deg-dist}}

For the first part of the corollary we show the result for $X^r_I$, the proof for $X_O^r$ being identical. We only need the result for $r=0$ and $r=1$ but the result holds for any $r$. Fix $r$. We use induction. The jump rate of $X_I$ if currently in state $i$ equals $\lambda^I_i=(i+\delta_I)c_I$. We start with $i=r$ meaning that there must be no jump between 0 and $t$:
$$
P(X_I^r(t)=r)=e^{-\lambda^I_rt}= e^{-(r+\delta_I)c_It},
$$
which agrees with the theorem. We now assume the expression for $P(X_I^r(t)=k)$ is as stated in the corollary for $k=1\dots ,i-1$ for all $t$. Then
\begin{align*}
P(X_I^r(t)&=i)  =\int_0^t P(X_I^r(s)=i-1)\lambda^I_{i-1}e^{-\lambda_i^I(t-s)}ds
\\
& = \frac{\delta_I+i-1)_{i-r}}{(i-1-r)!}c_Ie^{-(\delta_I+i)c_It} \int_0^t e^{(i-r)c_Is}(1-e^{-c_Is})^{i-1-r}ds
\\
&=\frac{\delta_I+i-1)_{i-r}}{(i-1-r)!}c_Ie^{-(\delta_I+i)c_It}  \sum_{j=0}^{i-1-r} \binom{i-1-r}{j} (-1)^j\frac{1}{c_I(i-r-j)}(e^{(i-r-j)c_It}-1)
\\
&=\frac{\delta_I+i-1)_{i-r}}{(i-1-r)!}\frac{e^{-(\delta_I+i)c_It}}{i-r} \left( \sum_{j=0}^{i-1-r} \binom{i-1-r}{j} (-1)^j e^{(i-r-j)c_It} +(-1)^{i-r} \right)
\\
&= \frac{\delta_I+i-1)_{i-r}}{(i-1-r)!} \frac{e^{-(\delta_I+i)c_It}}{i-r}  (e^{c_It}-1)^{i-r}
\\
&= \frac{\delta_I+i-1)_{i-r}}{(i-r)!}e^{-(\delta_I+r)c_It} (1-e^{-c_It})^{i-r},
\end{align*}
which proves the first part of the corollary. The third equality is obtained by expanding $(1-e^{-c_Is})^{i-1-r}$ to its binomial terms and integrating each term, and the fourth equality comes from summing the "-1" terms.

Now to the second part of the corollary. As before, we only show it for $X_I^r(T)$ (the proof for $X_O^r(T)$ being identical). If currently in state $k$, the rate at which $X_I$ gives birth equals $\lambda_k^I=(\delta_I+k)c_I$ and the rate at which the process stops equals $1$ ($T\sim Exp(1)$). The probability for a birth before a stop is hence $\lambda_k^I/(\lambda_k^I+1)$. Since the process is Markovian we hence have
$$
P(X_I^r(T)=i)=\left( \prod_{k=r}^{i-1}\frac{\lambda_k^I}{\lambda_k^I+1}\right)  \frac{1}{\lambda_i^I+1},
$$
where the last factor comes from requiring that the process stops when in state $k$. Simple manipulation of this expression, using that $\lambda_k^I=(\delta_I+k)c_I$, gives the desired result.

Now to the final part of the corollary. We know that the starting configuration of $X$ is either $(0,1)$ and $(1,0)$ with probabilities $\alpha/(\alpha+\gamma)$ and $\gamma/(\alpha+\gamma)$, so by conditioning on the starting configuration and the stopping time $T\sim Exp(1)$, and using that $X_I(t)$ and $X_O(t)$ are independent given the starting configuration, we get
\begin{align*}
p_{ij} =P(X(T)=(i,j))
&= \frac{\alpha}{\alpha+\gamma} \int_0^\infty P(X^0_I(t)=i)P(X_O^1(t)=j)e^{-t}dt 
\\
&\qquad + \frac{\gamma}{\alpha+\gamma} \int_0^\infty P(X^1_I(t)=i)P(X_O^0(t)=j)e^{-t}dt .
\end{align*}
From the first part of the corollory it follows that the first integral above equals
$$
\frac{(\delta_I+i-1)_i (\delta_O+j-1)_{j-1}}{i!(j-1)!} \int_0^\infty e^{-\delta_Ic_It}(1-e^{-c_It})^i e^{-(\delta_O+1)c_Ot}(1-e^{-c_Ot})^{j-1} dt.
$$
By expanding both $(1-e^{-c_It})^i$ and $(1-e^{-c_Ot})^{j-1}$, and integrating term by term, the first sum of the corollary is obtained. The second term is treated identically.

\subsection{Proof of Corollary \ref{cor-Bollob}}

We have that $p_i=P(X_I(T)=i)=\alpha/(\alpha+\gamma)P(X_I^0(T)=i) + \gamma/(\alpha+\gamma)P(X_I^1(T)=i) $, and these probabilities are given in Corollary \ref{cor-lim-deg-dist}: 
$$
P(X_I^r(T)=i)= \frac{1}{(\delta_I+\frac{1}{c_I} + i)} \frac{(\delta_I+i-1)_{i-r}}{(\delta_I+\frac{1}{c_I} + i-1)_{i-r}}  .
$$
It is easy to show that, as $i\to\infty$, $(a+i)_{i-r}/(b+i)_{i-r} \sim i^{a-b}$. This implies that 
$$
P(X_I^r(T)=i)\sim i^{-(1+1/c_I)},
$$
and since the same expression holds for $r=0$ and $r=1$ it also holds for $p_i$. The proof for $q_j$ is obtained similarly. We now fix $j$ and study $p_{ij}=P(X(T)=(i,j))$ for large $i$. We start by considering $j=0$, and hence look at $p_{i0}$ for large $i$. For $j=0$, the process must start in state $(1,0)$ and the out-degree birth process must never have a birth, so
\begin{align*}
p_{i0} &=\frac{\gamma}{\alpha+\gamma} \left( \prod_{k=1}^{i-1}\frac{\lambda_k^I}{\lambda_k^I+\lambda_0^O+1} \right) \frac{1}{\lambda_i^I+\lambda_0^O+1}
\\
& = \frac{(\delta_I+i-1)_{i-1}}{(\delta_I+\frac{1}{c_I} + \frac{c_O}{c_I}\delta_O + i-1)_{i-1}} \sim i^{-(1+ \frac{1}{c_I} + \frac{c_O}{c_I}\delta_O )},
\end{align*}
where the second equality is obtained by plugging in the expressions for $\lambda_k^I$ and $\lambda_0^O$, and the asymptotic size follows from the above stated property of ($a+i)_{i-r}/(b+i)_{i-r}$.

For another fixed $j\ge 1$ we now show that it is of the same order. In order to reach $(i,j)$ where $i$ is large and $j$ is fixed (and hence small in comparison with $i$) the $j$ births of the out-degree process can occur for different values of the in-degree process. However, since the birth rate for the in-degree process increases with $i$, the probability that any out-degree birth takes place when the in-degree is greater than $n$ is $o(1/n)$. So, the more likely paths ending in $(i,j)$ (where $j$ is fixed and small and $i$ is large) are where the out-degree births take place when the in-degree is small. We  now compute the probability for one such path, and since there are only a fixed number of such paths, and all these paths have probability of the same order, we conclude that the over-all probability $p_{ij}$ is of the same order as the probability of one such (likely) path.

We assume that the initial configuration is $(0,1)$ (the other case is done equivalently). We now compute the probability of the path, namely first having all the $j$ out-degree births, followed by all the indegree births:
\begin{align*}
P&((0,1)\to (0,2)\dots (0,j)\to (1,j)\to (2,j)\to \dots (i,j)=
\\
&= \prod_{k=1}^{j-1}\frac{\lambda_k^O}{\lambda_0^I+\lambda_k^O+1}  \prod_{k=0}^{i-1} \frac{\lambda_k^I}{\lambda_k^I+\lambda_j^O+1} \frac{1}{\lambda_i^I+\lambda_j^O+1}.
\end{align*} 
For fixed $j$, the first product is a constant, and the second product and the last factor are of order $i^{-(1+ \frac{1}{c_I} + \frac{c_O}{c_I}\delta_O )}$, which is shown in the same way as was done for $p_{i0}$. So, since $j$ is fixed, together with the fact that $j$ births of the out-degree process happen when the in-degree is small with large probability, proves the last statement of the theorem.

\subsection{Proof of Theorem \ref{Theor-tail}}

The exact expression for $p_{ij}$ is given in Equation (\ref{p_ijexact}). It is easy to show that for large $i=n$, $P(X_I^1(t)=n)\sim P(X_I^0(t)=n)$, and for $j=\floor{sn^r}$, $P(X_O^1(t)=\floor{sn^r})\sim P(X_O^0(t)=\floor{sn^r})$, which implies that whether starting in state $(0,1)$ or $(1,0)$ has no effect on the tail behavior. Further, it is known and easily proven that $(a+n)_n/n!\sim n^a$. From these two observations, together with the exact expression (\ref{p_ijexact}), we have
\begin{equation}\label{first_app}
p_{n,\floor{sn^r}}\sim n^{\delta_I-1} (sn^r)^{\delta_O-1} \int_0^\infty e^{-(c_I\delta_I+c_O\delta_O+1)t} \left(1-e^{-c_It}\right)^{n} \left(1-e^{-c_Ot}\right)^{sn^r} dt.
\end{equation}

We now analyse the integral in (\ref{first_app}) writing $b=c_I\delta_I+c_O\delta_O+1$ using Laplace's method. The integral in (\ref{first_app}) can be written as
$$
\int_0^\infty e^{-nf_n(t)}dt,\text{ where }f_n(t)=(b/n)t-\log(1-e^{-c_It}) - sn^{r-1}\log(1-e^{-c_It}).
$$
Clearly, $f_n(t)>0$ for $t>0$, and as $t\to 0$ or $t\to\infty$, $f_n(t)\to +\infty$. From this it follows that the main contribution to the integral comes for $t$-values close to $t_n^{(min)}$ for which $f_n(t)$ is minimized. We hence Taylor expand $f_n(t)$ around $t_{n}^{(min)}$: 
\begin{align*}
f_n(t) &= f_n (t_n^{(min)}) +(t-t_n^{(min)}) f'_n(t_n^{(min)}) +\frac{(t-t_n^{(min)})^2}{2} f''_n(t_n^{(min)}) +o( (t-t_n^{(min)})^2)
\\
&= f_n (t_n^{(min)}) + 0  +\frac{(t-t_n^{(min)})^2}{2} f''_n(t_n^{(min)}) +o( (t-t_n^{(min)})^2).
\end{align*}
By differentiating $f_n(t)$ it follows that, for large $n$, $t_n^{(min)}=\max (1/c_I, r/c_O)\log n + o(\log n)$ and $f_n^{(min)}:=f_n (t_n^{(min)})= \frac{\log n}{n} b\max (1/c_I, r/c_O)+o(\log(n)/n)$. Further, $f''_n(t)=c_I^2 e^{-c_It}+n^{r-1}c_O^2 e^{-c_Ot}$, and $f''_n(t_n^{(min)})\sim c_I^2n^{-c_I\max } + c_O^2  n^{r-1-c_O\max} \sim n^{-1}$, where $\max := \max (1/c_I, r/c_O)$. Going back to the intergral we hence have
\begin{equation} \label{int_approx}
\int_0^\infty e^{-nf_n(t)}dt \sim e^{-nf_n^{(min)}} \int_0^\infty e^{-n (t-t_n^{(min)})^2f''_n(t_n^{(min)})/2}dt\sim e^{-nf_n^{(min)}},
\end{equation}
the last step is by identifying the normal density. We have the following approximation for Equation (\ref{first_app}):
$$
p_{n,\floor{sn^r}}\sim n^{\delta_I-1} (sn^r)^{\delta_O-1} e^{-nf_n^{(min)}} \sim n^{\delta_I-1} n^{r(\delta_O-1)} n^{-(c_I\delta_I+c_O\delta_O+1)\max (1/c_I, r/c_O)},
$$
and this is the first statement of the theorem. 

As for the second statement, we have that $p_{n,\floor{n^r}}\sim n^{\delta_I-1+r(\delta_O-1)-b \max (1/c_I, r/c_O)}$. The $r$-value which maximizes this is hence the same $r$ which maximizes $g(r):=r(\delta_O-1)-b \max (1/c_I, r/c_O)$. For $r< c_O/c_I$, we have $g'(r)=\delta_O-1$ and for $r>c_O/c_I$, $g'(r)= -1-c_I\delta_I/c_O-1/c_O<0$, using that $b=c_I\delta_I + c_O\delta_O+1$. So, if $\delta_O<1$ the maximum is obtained at $r=0$, and if $\delta_O>1$, $g(r)$ is maximized for $r=c_O/c_I$. Finally, if $\delta_O=1$, then $g(r)$ is constant up to $r=c_O/c_I$ after which has negative derivative, so the maximum is obtained for any $r\in [0,c_O/c_I]$, which completes the proof.

\subsection{Proof of Corollary \ref{cor-Bol-Rio}}

We start by specifying the original PA-model for an undirected network once again. Suppose that at time $k=0$ the network consists of one single node without any edge. For each $k=1,2,\dots , n$, one node with $m$ edges is added, and each of the edges is attached, independently, to a node $u$ with degree $i$ with probability proportional to $i$. We now prove the result previously derived using other methods by Dorogovtsev et al \cite{DMS00}. After step $n$, let $p_i^{(n)}=P(\text{random node has degree }i)$. We want to prove that
\begin{equation}
p_i^{(n)} \to p_i:=\frac{2m(m+1)}{i(i+1)(i+2)}, i=m, m+1, \dots \label{p_i-un}
\end{equation}
As $n$ tends to infinity we approximate the PA model by a continuous-time process by speeding up time, similarly to what was done with the DPA model. For each $n$ we speed up time by letting new nodes enter after a time step $1/n$ rather than after one unit of time. So, for each $n$ we define the whole preferential attachment process up to step $n$: $X(1), ..., X(n)$ by $X^{(n)}(t),\ 0\le t\le 1$, where $X^{(n)}(t):= X( \floor{ tn} )$, where as before $\floor{ tn}$ denotes the integer part of $tn$.

More specifically, let's consider a node selected randomly among the $n+1$ nodes at time $1$ (after step $n$ in the original time scale), and analyse the distribution of its degree at time $1$. This degree will depend on when it was born (=entered the network), and since one node enters each time unit in the original time scale, the birth time is $U[0,\ 1]$ on the new time scale.

First we condition on the birth time $s$. Let $D^{(n)}(t),\ s\le t\le 1$ denote the degree of this node from birth until time 1 (when there are in total $n+1$ nodes). Since a node has degree $m$ when it enters the network (expect the first node which is selected with probability $1/(n+1)$ tending to 0) we have $D^{(n)}(s)=m$.  We hence seek the limit $\lim_n P(D^{(n)}(1)=i| D^{(n)}(s)=m)=: p_i(s)$.

At any time step, the total number of nodes and the total number of edges is non-random. At time $t$ ($\floor{ tn}$ in the original time scale) there are $\floor{ tn}+1$ nodes and $m\cdot \floor{ tn}$ edges. Since each edge is connected to two nodes, the sum of all degrees then equals $2m\floor{ tn}$. 

The process $D^{(n)}(t),\ s\le t\le 1$ is a pure birth chain with possible jumps at the increments $1/n$, so in the limit it converges to a continuous time time-inhomogeneous Markov birth process. Since $m$ edges are added when a new node is added, the degree could of course increase by more than 1 at a given time instant, but the probability of increasing by more than one is $O(1/n^2)$ and hence neglected. However, the probability to increase by 1 is $m$ times the probability that a specific edge connects to our node, plus terms of smaller order. We now compute the jump probability/intensity, which follows straightforward from the model definition. As before, let $\Delta t=1/n$. We have
\begin{equation}
\lambda^{(n)}_i(t)=P(D^{(n)}(t+\Delta t)=i+1|D^{(n)}(t)=i)/\Delta t \approx m \frac{i}{2m\floor{ tn}}\approx \frac{i/2}{t}  =: \lambda_i/t,
\end{equation}
where $\lambda_i=i/2$. Our limiting birth process is hence born at a uniform time $s$, and then evolves with birth rate $\lambda_i/t$ up until $t=1$. Similarly to the DPA model, this ditribution is equal to the distribution of a time-homogeneous linear birth process $Z(t)$ (a Yule process), born at time 0 in state $m$, having birth rate $\lambda_i$, and being stopped after and $Exp(1)$ time $T$. The distribution of a Yule process with birth rate $\lambda_i=i/2$, starting in state $m$ and observed at time $t$, has distribution
$$
P(Z(t)=k|Z(0)=m)=\binom{m-1}{k-1} e^{-mt/2}(1-e^{-t/2})^{k-1},
$$
e.g.\ de La Fortelle \cite{F06}.

Furthermore, if the birth process has intensities $\lambda_k=k/2$ and $T\sim Exp(1)$, then
\begin{align*}
P(Z(T)=i|Z(0)=m) &= \frac{\lambda_m}{\lambda_m+1}\cdot \dots \cdot \frac{\lambda_{i-1}}{\lambda_{i-1}+1}\frac{1}{\lambda_i+1}
\\
& =\frac{m/2}{m/2+1} \cdot \dots \cdot \frac{(i-1)/2}{(i-1)/2+1}\frac{1}{i/2+1} 
= \frac{2m(m+1)}{i(i+1)(i+2)},
\end{align*}
where the last equality is shown by induction.

\subsection{Proof of Theorem \ref{ThGDPA}}

The proof of Theorem \ref{ThGDPA} is more or less identical to the proof of Theorem \ref{mainth}. Now the limiting process can have births of each type, but also a simultaneous birth of both types (corresponding to the event that a new node attaches a directed edge, in either direction, to an existing node and the edge is reciprocated, which happens with probability $\rho$). This makes the limiting distribution $P(Y(T)=(i,j))$ harder to derive in that there are more paths to reach the state $(i,j)$, but the proof goes through in the same way.

\section{Conclusions and discussion}
In the current paper we analysed the directed preferential attachment (DPA) model using a new approach and showed that the limiting degree distribution can be characterized by two \emph{independent} (except starting either in state $(0,1)$ or $(1,0)$) birth processes  that are observed at a common $Exp(1)$ random time point, thus creating dependence. Beside shedding more light to the structure of the limiting degree distribution, this method also allowed for analyses of the bivariate tail probabilities (where both in and out-degree were assumed large), cf.\ Theorem \ref{Theor-tail}.

We also extended the DPA model to a more general model, the GDPA model, where new directed edges select nodes to attach to, in a way that may depend on both types of degrees, and, perhaps even more important, that the network may have a substantial fraction of edges being bi-directed (or equivalently undirected). The limiting degree distribution of this process was derived using similar methods. Also here, the distribution may be described by a bivariate birth process, but now the birth processes no longer evolve independently, and both components may have births at the same point in time. As before, the limiting degree distribution is the state of this process observed at an $Exp(1)$ random time point. The limiting tail probabilities for the GDPA model may perhaps also be derived using similar methods as in the DPA model, but this remains an open problem.

Beside deriving the tail distribution for the GDPA model it would be interesting to study other generalizations of the DPA network model. For instance, many empirical networks exhibit clustering, meaning that triangles are more common than for instance in the DPA model. It would be of interest to study related models that allow for such clustering to be present in the network, and to see if a similar alternative construction of the limiting tail behaviour may be obtained.

\section*{Acknowledgements}

I thank Sid Resnick and Svante Janson for valuable feedback during the early stage of this work, and Anders Martin-Löf for help with Laplace's method. I am grateful to the Swedish Research Council (grant 2015-05015) for financial support.

\end{document}